\documentclass[11pt, a4paper,leqno]{amsart}
\usepackage{amsmath,amsthm,amscd,amssymb,amsfonts, amsbsy}
\usepackage[english]{babel}
\usepackage{latexsym}
\usepackage{txfonts}
\usepackage{exscale}
\usepackage{enumitem}
\usepackage{bbm}
\usepackage{enumitem}
\usepackage{mathabx}
\usepackage{cite}

\usepackage[colorlinks,citecolor=red,hypertexnames=false]{hyperref}

\usepackage{color}

\calclayout
\allowdisplaybreaks

\theoremstyle{plain}
\newtheorem{theorem}{Theorem}
\newtheorem{lemma}[theorem]{Lemma}

\newtheorem{proposition}[theorem]{Proposition}

\theoremstyle{definition}

\theoremstyle{remark}
\newtheorem{remark}[theorem]{Remark}

\newcommand{\R}{\mathbb{R}}

\newcommand{\N}{\mathbb{N}}

\DeclareMathOperator{\Rp}{Re}
\DeclareMathOperator{\Ip}{Im}
\DeclareMathOperator{\sign}{sign}

\author[S.E. Boutiah]{Sallah Eddine Boutiah}
\author[L. Caso]{Loredana Caso}
\author[F. Gregorio]{Federica Gregorio}
\author[C. Tacelli]{Cristian Tacelli}

\address{Department of Mathematics, University Ferhat Abbas Setif-1, Setif 19000, Algeria\\ Laboratoire Des Math\'ematiques Appliqu\'ees-LaMA } 
\email{boutiah@univ-setif.dz} 

\address{Dipartimento di Matematica, Universita' degli Studi di Salerno, Via Giovanni Paolo II, 132, 84084 Fisciano (Sa), Italy}
\email{lorcaso@unisa.it}

\address{Dipartimento di  Matematica , Universita' degli Studi di Salerno, Via Giovanni Paolo II, 132, 84084 Fisciano (Sa), Italy}
\email{fgregorio@unisa.it} 

\address{Dipartimento di Matematica, Universita' degli Studi di Salerno, Via Giovanni Paolo II, 132, 84084 Fisciano (Sa), Italy}
\email{ctacelli@unisa.it}

\thanks{The second, the third and the fourth authors are members of Gruppo Nazionale per l'Analisi Matematica, la Probabilit\'a e le loro Applicazioni (GNAMPA) of the Istituto Nazionale di Alta Matematica (INdAM).  The first author is partially supported by The General Directorate of Scientific Research and Technological Development  (DGRSDT).  The authors wish to thank Abdelaziz Rhandi for several interesting discussions and useful comments.} 

\makeindex

\begin{document}

\title[Second-order operators with polynomially growing coefficients]{Some results on second-order elliptic operators with polynomially growing coefficients in $L^p$-spaces}

\maketitle

\begin{abstract}
In this paper we study minimal realizations in $L^p(\R^N)$ of the second order elliptic operator

\begin{equation*}\label{main oper}
{ A_{b,c}} := (1+|x|^\alpha)\Delta  + b|x|^{\alpha-2}x\cdot\nabla  - c |x|^{\alpha-2} - |x|^{\beta} , \quad x \in \mathbb{R}^N, 
\end{equation*}
where $N\geq3$, $\alpha\in[0,2)$, $\beta >0$, and $b, c$ are real numbers.  We use quadratic form methods to prove that $\left(A_{b,c},C_c^\infty\left(\R^N\setminus \{0\}\right)\right)$ 
admits an extension that generates     an  analytic   $C_0-$semigroup for  all $p\in(1,\infty)$. Moreover,    we give  conditions on the coefficients under which  
this extension is precisely the closure of $\left(A_{b,c},C_c^\infty\left(\R^N\setminus \{0\}\right)\right)$.

\end{abstract}

\section{Introduction}
\noindent  Second order elliptic operators with bounded coefficients have been widely studied in the
literature and nowadays they are
well understood, see the recent monograph \cite{LR20}. In recent years there have been substantial developments in the theory of second order elliptic operators  with unbounded coefficients. The latters arise as a model in many field of science, especially in stochastic analysis and mathematical finance, where stochastic models lead to equations with unbounded coefficients, e.g., the well known Black-Scholes equation.
 The developments regarding operators with unbounded coefficients are well documented, see  for example  \cite{can-rhan-tac1,boutiah-et-al2,can-rhan-tac2,can-tac1,dur-man-tac1,for-gre-rha,for-lor,kunze-luc-rha2,kunze-luc-rha1,lor-rhan,met-oka-sob-spi,met-spi2,met-spi3,met-spi4,met-spi-tac,met-spi-tac2,can-gre-rhan-tac,can-pap-tar,ber-lor} and the references therein.          

The present paper is devoted to the study of the following  elliptic operator

\begin{equation*}\label{main oper}
A_{b,c} u(x) = (1+|x|^\alpha)\Delta u(x) + b|x|^{\alpha-2}x\cdot\nabla u(x) -  c |x|^{\alpha-2}u(x) - |x|^{\beta} u(x), \quad x \in \mathbb{R}^N,
\end{equation*}
where $N \geq3$, $\alpha \in [0,2)$, $\beta> 0$ and $b, c \in \mathbb{R}$.  

 The case $\alpha >2$ has been treated in \cite{boutiah-et-al1}. In particular, the authors studied the complete operator   under the assumptions $\alpha >2$, $\beta >\alpha-2$, $b\in \R$  and  $c=0$. The main result concerns   
the generation by   a suitable realization of  the operator in $L^p(\R^N)$ of a strongly continuous and analytic semigroup   for $1<p<\infty$. 
See also  \cite{BelhajAli} and \cite{for-gre-rha} for the case $\alpha\geq0$ and $\beta=-2$ and \cite{can-rhan-tac1} for the case $b=c=0$.  

The operator $A_{b,c}$ without potential term $ c |x|^{\alpha-2} + |x|^{\beta}$ is studied in \cite{for-lor} for  $\alpha \in (1,2]$ and in \cite{met-spi-tac} for $\alpha>2$. 
In \cite{for-lor}  the  generation of an analytic semigroup in $L^p(\mathbb{R}^N)$  is proved for  $1\le p\le \infty$. 
 Whereas for $\alpha>2$, in \cite{met-spi-tac},  the generation of an analytic semigroup in $L^p(\R^N)$ is stated for $p>\frac{N}{N-2-b}$.

The Schr\"odinger-type operators $(1+ |x|^\alpha)\Delta- |x|^{\beta}$ in $L^{p}(\mathbb{R}^{N})$ have been considered in  \cite{lor-rhan} and \cite{can-rhan-tac1}.
In \cite{lor-rhan} the generation of an analytic semigroup in $L^p(\R^N)$ for $\alpha \in [0,2]$ and $\beta\geq0$ is obtained. In \cite{can-rhan-tac1} the same result is obtained in the case $\alpha>2$ and $\beta>\alpha-2$.

The operator with only diffusion term was studied  in \cite{met-spi2} and \cite{met-spi-tac2}.   
 In \cite{met-spi2} the authors proved 
 that for $N\ge 3$, $p>\frac{N}{N-2}$ and $\alpha>2$,  the maximal realization of the operator $(1+ |x|^\alpha)\Delta$ in $L^p(\mathbb{R}^N)$ generates an analytic semigroup  which is contractive if and only if $p\ge \frac{N+\alpha-2}{N-2}$. The case $\alpha\leq2$ has been treated in \cite{met-spi-tac2} and in this case no constraints on $p$ appear.

 In \cite{met-oka-sob-spi} the authors showed that the operator $\vert x \vert^{\alpha}\Delta + b\vert x \vert^{\alpha-2}x\cdot \nabla - c\vert x \vert^{\alpha-2}$ generates a strongly continuous semigroup in $L^p(\mathbb{R}^N)$ if  $s_1+\min\{0,2-\alpha\}<\frac{N}{p}<s_2+\max\{0,2-\alpha\}$, where $s_i$ are the roots of the equation $c+s(N-2+b-s)=0$. Here one can note that the case $\alpha=2$ is special.

We observe also that methods and results are twofold for $\alpha\leq2$ or  $\alpha > 2$. In particular, the generation results for $\alpha>2$ depend upon the dimension $N$; except for the case when the potential term has the form $|x|^\beta$ with $\beta>\alpha-2$.

Our aim in this paper is to complete the picture considering the complete operator
\[A_{b,c}=(1+|x|^\alpha)\Delta + b|x|^{\alpha-2}x\cdot\nabla  -  c |x|^{\alpha-2} - |x|^{\beta} \]
 in the case $\alpha<2$. While one expects to easily handle  the case $\alpha<2$,  when looking for more information on the generated semigroup and the domain of the operator, some difficulties arise and constraints on the coefficients appear. This is because, with the assumptions  $\alpha<2$ and $\beta>0$, the operator $A_{b,c}$ has unbounded coefficients at infinity and local singularities in the drift
and potential term.  In order to treat this singularity in  the origin, we initially study the operator in the set ${\mathcal D}_0:=C_c^\infty\left( \R^N\setminus\{0\} \right)$. 

We consider the minimal realization of $A_{b,c}$ in $L^p(\R^N)$: Our main focus is to establish if an extension of $(A_{b,c},{\mathcal D}_0)$ is the generator of a $C_0-$semigroup and when this extension is exactly its closure. We stress here that because of the singularities some conditions on $p$ and $N$ arise if one asks for essential $m-$accretiveness.   This is not surprising, 
for instance, the Laplace operator $-\Delta$ on ${\mathcal D}_0$ is essentially $m-$accretive if $p\leq \frac{N}{2}$.  The Schr\"odinger operator $-\Delta+c|x|^{-2}$ defined on ${\mathcal D}_0$ is essentially $m-$accretive in $L^{2}(\R^N)$ 
if and only if $c\geq 1-\frac{(N-2)^2}{4}$ and $N>4$ (see \cite{kalf}), the constant $ \frac{(N-2)^2}{4}-1$ being the best possible.  See also \cite{DR13} for other similar results. 
As regards the $L^p(\R^N)$ setting Okazawa shows that $-\Delta+c|x|^{-2}$ is essentially $m-$accretive if $c\geq c_0(N,p)$ where $c_0(N,p)$
depends on $N$ and $p$ (see \cite{Okazawa84}). 

  For this reason we treat our operator by means of two different methods.  In order to obtain generation results in $L^p(\R^N)$ without constraints on the coeffients, we first study the operator via form methods.
After proving some useful estimates we show that the introduced bilinear form is associated to a suitable extension  of the operator $\left(A_{b,c},\mathcal D_0\right)$ that generates a strongly continous analytic semigroup on $L^2(\mathbb{R}^N)$.  We then prove that when $c\geq0$ this semigroup is sub-Markovian and therefore  extrapolates to a family of consistent   $C_0-$semigroups on $L^p(\R^N)$ for $1\leq p<\infty$.  The case $c<0$, that means that the diffusion and the potential terms have the same sign, is more involved. Here an approximation technique is needed to obtain a family of consistent   $C_0-$semigroups on $L^p(\R^N), 1< p<\infty$.
The generators of the semigroups $(e^{t A_{b,c}})_{t\geq0}$ in $L^p(\R^N)$, denoted by $(A_{b,c},D_p(A_{b,c}))$, are extensions of $\left(A_{b,c},\mathcal D_0\right)$.

  On the other hand, in order to investigate when $(A_{b,c}, D_p(A_{b,c}))$ coincides with the closure of $(A_{b,c},{\mathcal D}_0)$,  we apply a perturbation theorem for linear $m$-accretive operators due to N. Okazawa, \cite{okazawa1996}. 
We consider the Schr\"odinger operator $A_{0,0}=(1+|x|^\alpha)\Delta-|x|^\beta$ already
studied in \cite{lor-rhan} and perturbe it with the potential $|x|^{\alpha-2}$.
Then, by means of a transformation, we are able to add also the drift term. It emerges that   for  $1<p<\frac{N-\alpha}{2-\alpha}$   the operator $(A_{b,c}, D_p(A_{b,c}))$   is the
generator of an analytic $C_0-$semigroup and it coincides with
the closure of  $(A_{b,c},{\mathcal D})$ where $\mathcal D:=C_c^\infty(\R^N)$. This in turn implies that $(A_{b,c}, D_p(A_{b,c}))$ is also the closure of $(A_{b,c},{\mathcal D_0})$.  
In the limit case $p=\frac{N-\alpha}{2-\alpha}$ an additional condition is required:
$
 \left(\frac{N}{p}-2+\alpha  \right)\left( \frac{N}{p'} -\alpha+b\right)+c>0$.
We finally
give more information on the domain $D_p(A_{b,c})$.

The following theorem summarizes the main results of the paper.

\begin{theorem} Suppose that $N\geq3$, $\alpha\in[0,2)$, $\beta>0$,  and $b, c\in\R$. 
\begin{itemize}
\item[(i)] For every $1<p<\infty$, $\alpha\neq0$, there exists $(A_{b,c}, D_p(A_{b,c}))$, an extension of $\mathcal ( A_{b,c},\mathcal D_0)$, that generates an  analytic    $C_0-$semigroup in $L^p(\R^N)$.
\item[(ii)] If $1<p<\frac{N-\alpha}{2-\alpha}$
or else $p=\frac{N-\alpha}{2-\alpha}$ and 
$ \left(\frac{N}{p}-2+\alpha  \right)\left( \frac{N}{p'} -\alpha+b\right)+c>0$, then $(A_{b,c}, D_p(A_{b,c}))$ coincides with the closure of $(A_{b,c},{\mathcal D})$ and it generates an analytic 
$C_0-$semigroup. Moreover,  in this case $\mathcal D_0$ is a core and

$$
 D_p(A_{b,c})  \subset\ \{u\in L^p(\R^N)\cap W^{2,p}_{loc}(\R^N\setminus\{0\}):
 (1+|x|^\alpha)D^2u, | |x|^{\alpha-1}\nabla u|, |x|^{\alpha-2}u,|x|^\beta u\in L^p(\R^N)\}.
$$
\end{itemize}

\end{theorem}

The paper is organized as follows. In Section \ref{opviaform}  the operator $A_{b,c}$ is studied via quadratic form methods. 
Section \ref{maximality} is concerned to the essential m-accretiveness of the operator $A_{b,c}$ and to further study of the domain $D_p(A_{b,c})$.

\section{Generation results in $L^p(\R^N)$}\label{opviaform}

The purpose of this section is to obtain generation results in $L^p(\R^N)$ for every $1<p<\infty$. To this aim we study the  operator
$A_{b,c}$
 via quadratic form methods.

Let $N\geq3$, $\alpha\in{(0,2)}$, $\beta>0$, $b,c\in\R$ and 
 
\begin{equation*}
A_{b,c} = (1+|x|^\alpha)\Delta  + b|x|^{\alpha-2}x\cdot\nabla  - c |x|^{\alpha-2} - |x|^{\beta} , \quad x \in \mathbb{R}^N.
\end{equation*}
We observe that if $u,v\in \mathcal{D}_0=  C_c^{\infty}(\R^N\setminus\{0\})$ then the following holds

\begin{align}\label{assocop}
-\int_{\R^N} &A_{b,c}u\overline v\,dx\\&\nonumber=\int_{\R^N} \left((1+|x|^\alpha)\nabla u\cdot \nabla \overline v+(\alpha-b)|x|^{\alpha-2}x\cdot \nabla u \overline v
        +c|x|^{\alpha-2}u\overline v+|x|^{\beta}u \overline v\right)\,dx.
\end{align}
Hence we define the bilinear form

\begin{align}\label{form}
& \mathfrak{a}(u,v)\\\nonumber&=\int_{\R^N}\left( (1+|x|^\alpha)\nabla u\cdot \nabla \overline v+(\alpha-b)|x|^{\alpha-2}x\cdot \nabla u \,\overline v
        +c|x|^{\alpha-2}u\overline v+|x|^{\beta}u \overline v+\lambda u\overline v\right)\,dx,\\\nonumber
D(\mathfrak{a})&=\{ u\in H^1(\R^N)\,:\,
    \left( 1+|x|^\alpha \right)^\frac{1}{2}\nabla u,
    \left( |x|^{\alpha-2} \right)^{\frac{1}{2}}u,
    \left( |x|^{\beta} \right)^{\frac{1}{2}}u\in L^2\left( \R^N \right)\},
\end{align}
where $\lambda$ is a suitable positive constant that will be chosen later.

Let us compute $\mathfrak{a}(u,u)$.

\begin{align*}
&\mathfrak{a}(u,u)\\&=\int_{\R^N} \left((1+|x|^\alpha)|\nabla u|^2+(\alpha-b)|x|^{\alpha-2}x\cdot \Rp (\overline u\nabla u )
        +(c|x|^{\alpha-2}+|x|^{\beta}+\lambda)|u|^2\right)\,dx,\\
&\quad +i \int_{\R^N}(\alpha-b)|x|^{\alpha-2}x\cdot \Ip (\overline u\nabla u )\,dx
        \end{align*}
and since $\Rp (\overline u\nabla u )=|u|\nabla |u|=\frac{1}{2}\nabla |u|^2$, integrating by parts we obtain

\begin{align*}
&\Rp \mathfrak{a}(u,u)\\&=\int_{\R^N}\left( (1+|x|^\alpha)|\nabla u|^2
        +\left[\left(c-\frac{\alpha-b}{2}(\alpha-2+N)  \right)|x|^{\alpha-2}+|x|^{\beta}+\lambda\right]|u|^2\right)\,dx.
\end{align*}

An application of the classical Hardy inequality, where we have denoted the Hardy constant by $c_0=\left( \frac{N-2}{2} \right)^2$, gives

\begin{align*}\label{eq:stimabeta}
\Rp \mathfrak{a}(u,u)&\geq \int_{\R^N} |x|^\alpha|\nabla u|^2\,dx\ +\frac{1}{2}\int_{\R^N}|\nabla u|^2\,dx \\\nonumber&\quad
        +\int_{\R^N}\left[\frac{c_0}{2|x|^{2}}+\left(c-1-\frac{\alpha-b}{2}(\alpha-2+N)  \right)|x|^{\alpha-2}+\lambda\right]|u|^2\,dx\\\nonumber&\quad+\int_{\R^N}|x|^{\alpha -2}|u|^2\,dx+\int_{\R^N}|x|^{\beta}|u|^2\,dx.
\end{align*}
 From now on we fix $\lambda>0$ in such a way that  for every $x\in\R^N$ the following quantity is positive
\begin{equation*}\label{pos}
\frac{c_0}{2|x|^{2}}+\left(c-1-\frac{\alpha-b}{2}(\alpha-2+N)  \right)|x|^{\alpha-2}+\lambda>0.\end{equation*}
This implies in particular that, for every $u\in D(\mathfrak{a})$, the following inequalities hold:

\begin{align}
\int_{\R^N} |\nabla u|^2\,dx&\leq 2\Rp \mathfrak{a}(u,u) \label{eq:stim0} \\
\int_{\R^N} (1+|x|^{\alpha})|\nabla u|^2\,dx&\leq 3\Rp \mathfrak{a}(u,u) \label{eq:stim1} \\
\label{eq:stima2}
\int_{\R^N} |x|^{\alpha-2}|u|^2\,dx&\leq \Rp \mathfrak{a}(u,u)\\
\label{eq:stima3}
\int_{\R^N} |x|^\beta |u|^2\,dx&\leq  \Rp \mathfrak{a}(u,u).
\end{align}

Moreover, $\mathfrak{a}$ is accretive on $D(\mathfrak{a})$, and hence we can
associate to $\mathfrak{a}$ the norm

\begin{align*}
\|u\|^2_\mathfrak{a}&=\int_{\R^N}\left( 1+|x|^\alpha \right)|\nabla u|^2\\\quad&+\left( \left(c-\frac{\alpha-b}{2}(\alpha-2+N)  \right)|x|^{\alpha-2}+|x|^\beta+\lambda+1\right)|u|^2 dx.
\end{align*}

We can easily show that $\mathfrak{a}$ is continuous taking into account that

\begin{align*}
\left|\int_{\R^N} |x|^{\alpha-2}x\cdot \nabla u\overline v\,dx\right|
&\leq \int_{\R^N} |x|^{\alpha-1}|\nabla u||v|\,dx\\&\leq  \left(\int_{\R^N} |x|^\alpha|\nabla u|^2\,dx\right)^{\frac{1}{2}}
\left( \int_{\R^N} |x|^{\alpha-2}|v|^2\,dx \right)^{\frac{1}{2}}.
\end{align*}

We want now to show that $\mathfrak{a}$ is closed. 

\begin{proposition}
The form $\mathfrak{a}$ is closed on $D(\mathfrak{a})$. 
\end{proposition}
{\sc Proof.}
Let $u_n\in D(\mathfrak{a})$ such that $\|u_n-u_m\|_\mathfrak{a}\to 0$ as $n,m\to \infty$.
From estimate \eqref{eq:stim0} it follows that $u_n$ is a Cauchy sequence in $H^1(\R^N)$ and then it converges to $u$ in  $H^1(\R^N)$ and $u$ belongs to $H^1(\R^N)$.

As a consequence of  \eqref{eq:stim1} we have that $v_n=(1+|x|^\alpha)^{\frac{1}{2}}\nabla u_n$ is a 
Cauchy sequence in $L^2(\R^N)$.
Since it converges to $(1+|x|^\alpha)^{\frac{1}{2}}\nabla u$ a.e., it also converges in $L^2(\R^N)$ and we have

 $$\int_{\R^N} (1+|x|^\alpha)|\nabla u_n-\nabla u|^2\,dx\to 0$$ as $n\to \infty$.
Estimates \eqref{eq:stima2} and \eqref{eq:stima3} give that 

$$\int_{\R^N} |x|^{\alpha-2}|u_n-u|^2\,dx\to 0$$
and 

$$\int_{\R^N} |x|^{\beta}|u_n-u|^2\,dx\to 0$$ as $n\to \infty$.
Then  $u\in D(\mathfrak{a})$ and $\lim_{n\to \infty}\|u_n-u\|_\mathfrak{a}=0$. This implies the closedness of $\mathfrak{a}$.

\qed

The next proposition shows that $\mathcal D_0$ is a core for $\mathfrak{a}$ in $L^2(\R^N)$.   
\begin{proposition}\label{pr:core-for-a}
${\mathcal D}_0$ is a core for $\mathfrak{a}$.
\end{proposition}
{\sc Proof.}
Let us prove that $H^1_0(\R^N \setminus \{0\} )$ is dense in $D(\mathfrak{a})$ with respect to
the form norm. 
Take $u\in D(\mathfrak{a})$ and consider $u_n=u\varphi _n$,
where $\varphi_n\in \mathcal D_0$ is such that 

\begin{equation*}  
\left\{
\begin{array}{ll}
\varphi_n=0  \  \text{in} \   B({\frac1n})\cup B^c({2n}), \\
\varphi_n=1 \  \text{in} \  B(n)\setminus B({\frac2n)}, \\
0\leq \varphi_n\leq 1,\\
|\nabla \varphi_n(x)|\leq C\frac{1}{|x|}.
\end{array}\right.
\end{equation*}

Observe that $u_n\in H^1(\R^N)=H^1(\R^N\setminus\{0\})$, then since $N\geq3$, 
 by \cite[Remark 17]{Brezis}, $u_n\in H_0^1(\R^N\setminus\{0\})$. 

We show now that $u_n\to u$ with respect to the norm of $\mathfrak{a}$.
Indeed, $u_n|x|^{\beta/2}\to u|x|^{\beta/2}$ and $u_n|x|^{\frac{\alpha-2}{2}}\to u|x|^{\frac{\alpha-2}{2}}$ in $L^2(\R^N)$ by dominated convergence.

As regards the other term we observe that

\begin{eqnarray*}
(1+|x|^\alpha)^{\frac{1}{2}}|\nabla u_n-\nabla u| &=&
(1+|x|^\alpha)^{\frac{1}{2}}\left|(\varphi_n -1)\nabla u + u\nabla \varphi_n\right|\\
&\le & (1+|x|^\alpha)^{\frac{1}{2}}(1-\varphi_n )|\nabla u|+ (1+|x|^\alpha)^{\frac{1}{2}}|u||\nabla \varphi_n|.
\end{eqnarray*}
Therefore, it suffices to prove that the second term converges to 0 in $L^2(\R^N)$. Indeed it converges a.e. to $0$ since

\begin{align*}
 (1+|x|^\alpha)^{\frac{1}{2}}|u|  |\nabla \varphi_n|\leq
 C \frac{ (1+|x|^{\alpha})^{\frac{1}{2}} }{|x|}|u|\chi_{K_n}\leq
 |u|\left( C_1\frac{1}{|x|}+C_2 |x|^{\frac{\alpha-2}{2}}\right)\chi_{K_n}
\end{align*}
where $K_n=B({\frac2n)}\setminus B({\frac1n)}\cup B({2n})\setminus B(n)$. We also observe that 
$\int_{\R^N} \frac{|u|^2}{|x|^2}dx<\infty$ by Hardy's inequality and 
$\int_{\R^N} |u|^2|x|^{\alpha -2}dx<\infty$ since $u\in D(\mathfrak{a})$. Then 
$(1+|x|^\alpha)^{\frac{1}{2}}|u|  |\nabla \varphi_n|$ converges to $0$ in $L^2(\R^N)$ by dominated convergence.

\qed

As a consequence of the previous results we obtain the following generation theorem that we state without proof, c.f. \cite[Theorem 1.52]{ouhabaz}.  

\begin{theorem}\label{thmform}
The form $\mathfrak{a}$ defined in \eqref{form}  is densely defined, accretive, continuous and closed. Therefore, it is associated to a closed  operator $(-\widetilde{  A}_{b,c},D(\widetilde A_{b,c}))$ on $L^2(\R^N)$ defined by

\begin{align*}
D(\widetilde{A}_{b,c}):&=\{u\in D(\mathfrak{a})\,:\,\exists\,v\in L^2(\R^N)\,s.t.\,\,\mathfrak{a}(u,h)=\langle v,h\rangle ,\,\forall\,h\in D(\mathfrak{a})\}\\-\widetilde{ A}_{b,c}u:&=v.
\end{align*}
Moreover, $(\widetilde{  A}_{b,c}, D(\widetilde {A}_{b,c}))$ is the generator of a strongly continuous analytic contraction semigroup on $L^2(\R^N)$.
\end{theorem}

It is clear that $(\widetilde{ A}_{b,c}, D(\widetilde {A}_{b,c}))$ is an extension of $(A_{b,c}-\lambda ,\mathcal D_0)$. Indeed if $u\in \mathcal{D}_0$, by \eqref{assocop} it follows that 

$$\mathfrak{a}(u,v)=-\int_{\R^N} (A_{b,c}u-\lambda u)v\,dx$$ for every $v \in \mathcal{D}_0 $ and, since $\mathcal{D}_0$  is a core for $\mathfrak{a}$, it holds for every $v\in D(\mathfrak{a})$.
Hence $A_{b,c}u-\lambda u=\widetilde A_{b,c} u$.

We now study further properties of the semigroup $(e^{t\widetilde A_{b,c}})_{t\geq0}$. If our analytic contraction semigroup on $L^2$ satisfies the  $L^\infty$-contractivity property, then using the Riesz-Thorin interpolation theorem, we can extend each operator from $L^2(\R^N)\cap L^p(\R^N)$ to a contraction   operator on $L^p$ for $2\le p\leq\infty$.   

\begin{theorem}\label{th:extrapolates-cgeq0}
 Let $N\geq3$, $\alpha\in(0,2)$, $\beta>0$, $b\in\R$ and $c\geq0$.
 The semigroup $(e^{t\widetilde A_{b,c}})_{t\geq0}$ extrapolates to a
  family of consistent 
 contraction analytic $C_0-$semigroups on $L^p(\R^N)$ for every $p\in[2,\infty)$. The generators $(\widetilde A_{b,c},D_p(\widetilde {A}_{b,c}))$ are   extensions of $(A_{b,c}-\lambda,\mathcal D_0)$.
\end{theorem}

{\sc Proof.} 
We prove that the semigroup $(e^{t\widetilde A_{b,c}})_{t\ge 0}$ is sub-Markovian, i.e. $(e^{t\widetilde A_{b,c}})_{t\ge 0}$ is both positive and $L^\infty$-contractive.   To this purpose we first observe that the semigroup is real.  By \cite[Proposition 2.5]{ouhabaz} the semigroup is real if 

\begin{itemize}
\item[(a)] $u\in D(\mathfrak{a})\Longrightarrow\,\Rp u \in D(\mathfrak{a}),\ \mathfrak{a}(\Rp u,\Ip u)\in \R$.
\end{itemize}

Then let $u\in D(\mathfrak{a})$, it follows that,

$$\int_{\R^N} (1+|x|^\alpha) |\nabla\Rp u|^2\,dx\le \int_{\R^N} (1+|x|^\alpha)|\nabla u|^2\,dx<\infty$$
and hence $(1+|x|^\alpha)^\frac12\nabla\Rp u\in L^2(\R^N)$.
One can make similar computations to show that  

$(|x|^{\alpha-2})^\frac12\Rp u\in L^2(\R^N)$ and  $(|x|^\beta)^\frac12\Rp u\in L^2(\R^N)$. The condition $\mathfrak{a}(\Rp u,\Ip u)\in \R$ is straightforward, then $(e^{t\widetilde A_{b,c}})_{t\ge 0}$ is real. 

Now, by \cite[Theorem 2.6]{ouhabaz} the semigroup is positive if and only if 
\begin{itemize}
\item [(b)] the semigroup is real and for every $u\in D(\mathfrak{a})\cap L^2(\R^N)$ real-valued we have  $u^+\in D(\mathfrak{a})$ and  $\mathfrak{a}(u^+,u^-)\leq0$.
\end{itemize} 
It suffices then to notice that $\mathfrak{a}(  u^+, u^-)=0$.

Finally, we can apply \cite[Corollary 2.17]{ouhabaz} that states that the sub-Markovian property is equivalent to
\begin{itemize}
\item[(c)] for  $u\in D(\mathfrak{a})\cap L^2(\R^N)^+\ \Longrightarrow 1\wedge u\in D(\mathfrak{a}), \mathfrak{a}(1\wedge u,(u-1)^+)\geq0$.
\end{itemize}
Recall the notations $1\wedge u  := \inf\{1, u \}$ and $(u -1)^+ := \sup\{ u -1,0\}$. 
Then let $u\in D(\mathfrak{a})$ be real and non-negative. It is known that  $(1\wedge u ) \in H^1(\R^N)$  then, it follows that,
$$\int_{\R^N} (1+|x|^\alpha)|\nabla (1\wedge u )|^2\,dx\le \int_{\R^N} (1+|x|^\alpha)|\nabla u|^2\,dx<\infty$$
and hence $(1+|x|^\alpha)^\frac12|\nabla( 1\wedge u ) |\in L^2(\R^N)$. A similar argument  leads to $(|x|^{\alpha-2})^\frac12(1\wedge u ) \in L^2(\R^N)$ and  $(|x|^\beta)^\frac12(1\wedge u ) \in L^2(\R^N)$. 

Moreover one has

\begin{align*}\label{sign}
 \mathfrak{a}(1\wedge u,(u-1)^+)
        &=\int_{ u \geq1}(c|x|^{\alpha-2}+|x|^\beta+\lambda) (u-1)\,dx\geq0.
\end{align*}

Therefore one obtains $||e^{t\widetilde A_{b,c}}||_\infty\leq1$.
Taking into account the contractivity  on $L^2$ and using the Riesz-Thorin interpolation theorem 
we can extend each operator $e^{t\widetilde A_{b,c}}$ to an operator on $L^p(\R^N)$  for every $2\leq p\leq \infty$.
Moreover the family $( e^{t\widetilde A_{b,c}} )_{t\geq 0}$ defines a consistent family of $C_0-$semigroup of   contractions on $L^p(\R^N)$ for every $2\leq p<\infty$.
\qed
\\

Let us now focus on the case $c<0.$
Consider  the extrapolated semigroups when $c=0$, $(e^{t  \widetilde A_{b,0}})_{t\geq 0}$,
and consider the sequence of bounded potentials $W_n=\max \{-n,c|x|^{\alpha-2}\}$ with $c<0$. By perturbation of bounded operators we have that the sum $   \widetilde A_{b,0}-W_n$ generates a
$C_0-$semigroup on $L^p(\R^N)$. Moreover

$$0\leq e^{t(  \widetilde A_{b,0}-W_n)}\leq e^{t(  \widetilde A_{b,0}-W_{n+1})}.$$

Reasoning as in \cite{Arendt-et-al} one can prove the following result.
\begin{theorem} \label{th:extrapolates-cleq0}
Let $N\geq3,\alpha\in(0,2), \beta>0,b\in\R$ and $c<0$. 
The sequence of semigroups $e^{t(  \widetilde A_{b,0}-W_n)}$ converges to a contractive $C_0-$semigroup in $L^p(\R^N)$ 
and its generator $(\widetilde A_{b,c},D_p(\widetilde {A}_{b,c}))$ is an extension of $(A_{b,c}-\lambda,\mathcal D_0)$.
\end{theorem}

{\sc Proof.}
To prove the convergence of the semigroups we proceed as in \cite[Proposition 3.6]{Arendt-et-al}.
Since the sequence is monotone we have to prove that the semigroups $e^{t  (  \widetilde A_{b,0}-W_n)}$ are contractive.
For that we prove the dissipativity of $\widetilde A_{b,0}-W_n$ for all $n \in \mathbb{N}$.
We emphasize the dependence on $c$ setting $\mathfrak{a}=\mathfrak{a}_c$.
Let $u\in D_p(  \widetilde A_{b,0}-W_n)=D_p(  \widetilde A_{b,0})$. By \cite[Theorem 3.11]{ouhabaz}  one has
$|u|^{p/2}\in D(\mathfrak{a}_0)$ and 

\begin{align*}
\Rp(-  \widetilde A_{b,0}u+W_nu,|u|^{p-1}\sign u) &\geq \frac{2}{p} \mathfrak{a}_0(|u|^{p/2},|u|^{p/2})+\Rp(W_nu,|u|^{p-1}\sign u)\\
&=\frac{2}{p}\mathfrak{a}_{cp/2}(|u|^{p/2},|u|^{p/2})\geq 0.
\end{align*}

Hence, $\widetilde A_{b,0}-W_n$ is dissipative for all $n \in \mathbb{N}$ 
and  it follows that the semigroups $e^{t (  \widetilde A_{b,0}-W_n)}$ are contractive  and converge strongly in $L^p(\R^N)$ 
to a $C_0-$semigroup which we denote by $e^{t B}$.

Now we want to show that $(B,D(B))$ is an extension of $(A_{b,c}-\lambda,{\mathcal D}_0)$.
Let $u\in {\mathcal D}_0$. 
We have $u\in D_p(\widetilde A_{b,0})$. Moreover
there exists $n_0$ such that for every $n\geq n_0$ we have $W_n=W$ in the support of $u$.
Then
\[(  \widetilde A_{b,0}-W_n)u=(\widetilde A_{b,0}-W)u=\widetilde { A}_{b,c}u=v.\]
Now we observe that for all $n\geq n_0$
\[
e^{t  (  \widetilde A_{b,0}-W_n)}u-u=\int _0^t e^{s (  \widetilde A_{b,0}-W_n)} vds.
\]
 Taking the limit for $n\to \infty$ we have
\[
e^{t B}u-u=\int _0^t e^{s  B} vds
\]
then $u\in D(B)$ and $\widetilde { A}_{b,c}u= Bu.$
\qed

Therefore, we obtained  the following theorem.

\begin{theorem}\label{formgeneration} Let $N\geq3,\alpha\in(0,2),\beta>0, b,c\in \R$. 
There exists $ (A_{b,c},D_p(A_{b,c}))$, an extension of $( A_{b,c},\mathcal D_0)$, that generates an   analytic  $C_0-$semigroup in $L^p(\R^N)$ for any $1<p<\infty$. 
\end{theorem}
{\sc Proof.}
By Theorem \ref{th:extrapolates-cgeq0} and Theorem \ref{th:extrapolates-cleq0} we have that 
for every $c\in \R$,
$(e^{t\widetilde A_{b,c}})_{t\geq0}$ extrapolates to a
  family of consistent 
 contraction $C_0-$semigroups on $L^p(\R^N)$ for every $2\leq p<\infty $.
 By duality we have that 
$(e^{t\widetilde A^{*}_{b,c}})_{t\geq0}$ defines a
  family of consistent 
 contraction $C_0-$semigroups on $L^p(\R^N)$ for every $1<p\leq 2$.
 Now we observe that  the adjoint 
 $(  \widetilde A^*_{b,c},D(  {\widetilde A}_{b,c}^*))$ is an extension of $(\widetilde A_{b_1,c_1},{\mathcal {D}_0})$ where 

\[\widetilde A_{b_1,c_1}=(1+|x|^\alpha)\Delta+b_1|x|^{\alpha-2}x\cdot\nabla-c_1|x|^{\alpha-2}-|x|^\beta-\lambda\]
with $b_1=2\alpha-b$ and $c_1=(b-\alpha)(\alpha-2+N)$.

Then, rearranging the coefficients, we obtain that $(e^{ t \widetilde A_{b,c}})_{t\geq 0}$ is $L^{p}-$contractive for every $1< p<\infty$.
 
 Taking into account the analyticity of the semigroup $(e^{t{\widetilde A}_{b,c}})_{t\geq0}$ in $L^2(\R^N)$ and the  contractivity   in $L^p(\R^N)$ for every $1< p<\infty$, by \cite[Proposition 3.12]{ouhabaz}, the semigroups $(e^{t\widetilde A_{b,c}})_{t\geq0}$ are analytic in $L^p(\R^N)$ for every $1< p<\infty$.

Now, since $(\widetilde{  A}_{b,c}, D_p(\widetilde{A}_{b,c}))$ is an extension of $\left(A_{b,c}-\lambda,\mathcal D_0\right)$, 
we can state that an extension of $(A_{b,c},{\mathcal D}_0)$, denoted by   $ (A_{b,c},D_p({A}_{b,c}))$, generates
the   analytic  $C_0-$semigroup $e^{t A_{b,c}}=e^{\lambda t}e^{t\widetilde A_{b,c}}$ 
  in $L^p(\R^N)$ for $1<p<\infty$.

\qed

\begin{remark}
\begin{itemize}
\item   Observe that when $c\geq0$ thanks to the $L^\infty$ contractivity the semigroup is strongly continuous also in $L^1(\R^N)$. 
\item
The generation of an analytic semigroup in $L^p(\R^N)$, $p\in(1,\infty)$ when $\alpha\in[1,2)$, $\beta\geq0$, $N\geq3$ and $c=0$ is straightforward. Indeed, in this case the operator $F=b|x|^{\alpha-2}x\cdot\nabla$ is a small perturbation of the operator 
$A_{0,0}=(1+|x|^\alpha)\Delta-|x|^\beta$
studied in \cite{lor-rhan}. It is easy  to prove then that  the sum operator

\[A_{0,0}+F=(1+|x|^\alpha)\Delta+b|x|^{\alpha-2}x\cdot\nabla-|x|^\beta\]
with domain

\begin{equation}\label{doma0}D(A_{0,0})=\{u\in W^{2,p}(\R^N)\,:\,(1+|x|^\alpha)|D^2 u|,(1+|x|^\alpha)^\frac12 |\nabla u|,|x|^\beta u\in L^p(\R^N)\}\end{equation}
generates an analytic semigroup on $L^p(\R^N)$ for $1<p<\infty$.

Indeed, it suffices to notice that by \cite[Theorem 2.5 and Proposition 2.9]{lor-rhan} we know that the $L^p$ realisation of $A_{0,0}$ with domain \eqref{doma0} generates a strongly continuous analytic semigroup in $L^p(\R^N)$, $p\in(1,\infty)$. Since for $\alpha\in[1,2)$ the term $|x|^{\alpha-1}$ is controlled by $(1+|x|^\alpha)^\frac12$ one has that $D(A_{0,0})\subset D(F)$. Moreover, we can apply the estimate (2.3) of \cite[Proposition 2.3]{lor-rhan} to infer that there exist two positive constants $\varepsilon_0, C_\varepsilon$ such that

\begin{equation}\label{A bounded}
||Fu||_p\leq||(1+|x|^\alpha)^\frac12|\nabla u|||_p
\leq\varepsilon||(1+|x|^\alpha)\Delta u-|x|^\beta u||_p+C_\varepsilon||u||_p
\end{equation}
for any $\varepsilon\in(0,\varepsilon_0]$ and $u\in W^{2,p}(\R^N)$.
By \eqref{A bounded} we obtain that $F$ is $A_{0,0}$-bounded with $A_{0,0}$ bound as small as we want. Therefore we can apply \cite[III, Theorem 2.10]{engel-nagel} to claim that $(A_{0,0}+F, D(A_{0,0}))$ generates an analytic semigroup in $L^p(\R^N)$ for any $1<p<\infty$.
\end{itemize}
\end{remark}

\section{Essential m-accretivity}\label{maximality}

 This section is concerned with the question 
 when $(A_{b,c}, D_p(A_{b,c}))$ coincides with the closure of $(A_{b,c},{\mathcal D}_0)$. In order to obtain this result  one has to investigate on essential maximal accretivity of the operator. We first deal with the case $b=0$ and then, by means of a transformation, we add the drift term. 

\subsection{The case $b=0$}

In this subsection we consider the case $b=0$ and we treat the  operator as a sum of 
$A_{0,0}=(1+|x|^\alpha)\Delta-|x|^{\beta}$ and the potential $W=|x|^{\alpha-2}$
and show that $-A_{0,c}=-A_{0,0}+cW$ defined on  

\begin{equation*}\label{dom}
 D(A_{0,0})=\{u\in W^{2,p}(\R^N)\,:\,(1+|x|^\alpha)|D^2 u|,(1+|x|^\alpha)^\frac12 |\nabla u|,|x|^\beta u\in L^p(\R^N)\}\end{equation*}
is maximal accretive on $L^p(\R^N)$ for $1<p<\frac{N-\alpha}{2-\alpha}$ and every $c\in \R$. Moreover 
$C_c^\infty(\R^N)$ is a core for $(A_{0,c},{\mathcal D}_0)$.

As already recalled, in \cite{lor-rhan}  L. Lorenzi and A. Rhandi proved generation results for $A_{0,0}$ and that
 $C_c^{\infty}(\R^N)$ is a core for $(A_{0,0},  D(A_{0,0}))$.
In order to treat $-A_{0,0}+cW$ we use the 
following perturbation theorem due to N. Okazawa, see \cite[Theorem 1.7]{okazawa1996}.

\begin{theorem}\label{Okazawa 1.7}
Let $A$ and $B$ be linear $m-$accretive operators in $L^p(\R^N)$. Let $D$ be a core for $A$
and let $\{B_\varepsilon\}$ be the Yosida approximation of $B$.
Assume that
\begin{itemize}
\item[(i)] there are constants $a_1,a_2\ge 0$ and $k_1>0$ such that for all $u\in D$ and $\varepsilon>0$

\[
\Rp \langle Au,F(B_\varepsilon u) \rangle\geq k_1 \|B_\varepsilon u\|^2_p
-a_2\|u\|_p^2-a_1\|B_\varepsilon u\|_p\|u\|_p
\]
where $F$ is the duality map on $L^p(\R^N)$ to $L^{p\prime}(\R^N).$
\end{itemize}
Then $B$ is $A-$bounded with $A-$bound $k_1^{-1}$:

\[
 \|Bu\|_p\leq k_1^{-1}\|Au\|_p+k_0\|u\|,\quad u\in D(A)\subset D(B).
\]
Assume further that
\begin{itemize}
\item[(ii)] $\Rp \langle u,F(B_\varepsilon u)  \rangle\geq 0$, for all $u \in L^p(\R^N)$ and $\varepsilon>0$;
\item[(iii)] there is $k_2>0$ such that $A - k_2 B$ is accretive.   
\end{itemize}
Set $k = \min\{k_1,k_2 \}$. If $t >-k$ then $A + t B$ with domain $D(A+t B) = D(A)$ is m-accretive
and any core of $A$ is also a core for $A + t B$.
Furthermore, $A - k B$ is essentially m-accretive on $D(A)$. 
\end{theorem}

In our situation  the Yosida approximation of $W$ is given by
$W_\varepsilon(x)=\frac{W}{1+\varepsilon W}=\frac{1}{|x|^{2-\alpha}+\varepsilon}$ 
and the duality map is given by
$F(W_{\varepsilon})=W_{\varepsilon}^{p-1}\overline u|u|^{p-2}\|W_\varepsilon u\|_p^{2-p}.$ 

We will use a $L^p-$generalization of the classical Hardy inequality. If $1<p<\infty$ and  $\delta \ge 0$ then the following inequality holds

\begin{equation}\label{eq:gen-hardy-weight}
 \left( \frac{N-2+\delta}{p} \right)^2\int_{\R^N} \frac{|v|^p}{|x|^2}|x|^\delta dx\leq \int_{\R^N} |\nabla v|^2|v|^{p-2}|x|^{\delta}dx
\end{equation}
for every $v\in C_c^\infty\left( \R^N \right)$, cf. \cite[Lemma 2.4]{for-gre-rha}.  
In order to apply the above theorem we need the following estimates that we prove for the complete operator $A_{b,c}$.

\begin{lemma}\label{estimate-pot-eps}
Let $p>1$.  Assume that $N\geq3$, $0 \le \alpha < 2$, $\beta >0$. 

\begin{itemize}
\item[$\left(i \right)$] If $p<\frac{N-\alpha}{2-\alpha}$,
then for every  $u\in C_c^{\infty}(\R^N)$ and for every $k\in \R$ there exists $a_1\geq 0$ such that 

\begin{equation}\label{estim 1}
\Rp \langle -A_{b,c}u,W_{\varepsilon}^{p-1} u|u|^{p-2} \rangle\geq k \|W_\varepsilon u\|^p_p
-a_1\|W_\varepsilon u\|_p^{p-1}\|u\|_p.
\end{equation}
\item[$\left(ii \right)$] If $p=\frac{N-\alpha}{2-\alpha}$, 
then for every  $u\in C_c^{\infty}(\R^N)$ we have

\begin{equation}\label{estim 2}
\Rp \langle -A_{b,c}u,W_{\varepsilon}^{p-1} u|u|^{p-2} \rangle\geq k_1 \|W_\varepsilon u\|^p_p
-a_1\|W_\varepsilon u\|_p^{p-1}\|u\|_p
\end{equation}
where 
$k_1=k_1(b,c)=\left( \alpha-2+\frac{N}{p} \right)\left( -\alpha+b+\frac{N}{p'} \right)+c$ and $a_1\geq 0.$
\end{itemize}
\end{lemma}
{\sc Proof.}
Let $u\in C_c^{\infty}(\R^N)$.
Observe that,  thanks to the condition $p\leq \frac{N-\alpha}{2-\alpha}$,  the term  $|x|^{\alpha-2}$ is locally in $L^p(\R^N)$. Hence, when multiplied by a compactly supported function, the resulting function belongs to $L^p(\R^N)$. 

Multiplying ${A_{b,c}}u$ by $W_{\varepsilon}^{p-1}\overline u|u|^{p-2}$ and integrating over $\R^N$ we obtain the following.
We note here that the integration by part in the singular case $1<p<2$ is allowed thanks to \cite{met-spi}.

\begin{align*}
&-\int_{\R^N} {A_{b,c}}u W_{\varepsilon}^{p-1}\overline u|u|^{p-2} dx
\\&\quad =-\int_{\R^N} a(x)\Delta u W_{\varepsilon}^{p-1}\overline u|u|^{p-2}
 +b|x|^{\alpha-2}x\cdot \nabla u W_{\varepsilon}^{p-1}\overline u|u|^{p-2}dx\\
 &\qquad
    +\int_{\R^N}c|x|^{\alpha-2}W_{\varepsilon}^{p-1}|u|^{p}+|x|^{\beta}W_{\varepsilon}^{p-1}|u|^{p}dx, 
\end{align*}
where $a(x) = (1+|x|^\alpha)$. Then

\begin{align*}
&-\int_{\R^N} a(x)\Delta u W_{\varepsilon}^{p-1}\overline u|u|^{p-2}dx=
\int_{\R^N} \nabla u\cdot \nabla \left ( a(x) W_{\varepsilon}^{p-1}\overline u|u|^{p-2}\right )dx\\\nonumber
&\quad =\int_{\R^N}\left[ \nabla u \cdot \nabla \left( \overline u |u|^{p-2} \right)a(x) W_{\varepsilon}^{p-1}
	+a(x)\nabla u\cdot \nabla W_{\varepsilon}^{p-1}\overline u|u|^{p-2}
	\right.\\\nonumber&\qquad\left.+\alpha |x| ^{\alpha-2}x\nabla u W_{\varepsilon}^{p-1}\overline u|u|^{p-2}\right]\,dx.
\end{align*}

Recalling that

$$\Rp\left( \nabla u \cdot \nabla \left( \overline u |u|^{p-2} \right)\right)
=(p-1)|u|^{p-4}|\Rp\left(\overline u\nabla u\right)|^2
+|u|^{p-4}|\Ip\left(\overline u\nabla u\right)|^2, $$
and

 $$\Rp \left( \overline u\nabla u \right)=|u|\nabla |u|, $$ 
we have 

\begin{align*}
&\Rp\left(-\int_{\R^N} a(x)\Delta u W_{\varepsilon}^{p-1}\overline u|u|^{p-2}
 +b|x|^{\alpha-2}x\cdot \nabla u W_{\varepsilon}^{p-1}\overline u|u|^{p-2}
dx\right) \\
&\qquad =(p-1)\int_{\R^N} a(x) W_{\varepsilon}^{p-1} |u|^{p-4}|\Rp\left(\overline u\nabla u\right)|^2dx\\
&\qquad 
        +\int_{\R^N} a(x) W_{\varepsilon}^{p-1}|u|^{p-4}|\Ip\left(\overline u\nabla u\right)|^2dx\\
&\qquad  +\frac{1}{p} \int_{\R^N} a(x) \nabla W_{\varepsilon}^{p-1}\cdot\nabla |u|^p
	+(\alpha -b) W_{\varepsilon}^{p-1}|x| ^{\alpha-2}x \cdot \nabla |u|^p dx. 
\end{align*}
Then,

\begin{align*}
&\Rp\left(-\int_{\R^N} A_{b,c} u W_{\varepsilon}^{p-1}\overline u|u|^{p-2}dx\right)\\
&\quad \geq (p-1)\int_{\R^N} (1+|x|^\alpha) W_{\varepsilon}^{p-1} |u|^{p-2}|\nabla |u||^2dx\\
&\qquad + \int_{\R^N}\left[ -\frac{1}{p}{\rm div}\left[ a(x) \nabla W_{\varepsilon}^{p-1}
	+(\alpha-b) W_{\varepsilon}^{p-1}|x| ^{\alpha-2}x \right] |u|^p\right.\\&\qquad\left.+
	c|x|^{\alpha-2}W_{\varepsilon}^{p-1}|u|^{p}+|x|^\beta W_\varepsilon^{p-1}|u|^p \right]dx\\
&\quad = (p-1)\int_{\R^N} (1+|x|^\alpha) W_{\varepsilon}^{p-1} |u|^{p-2}|\nabla |u||^2dx\\
&\qquad  -\frac{1}{p}\int_{\R^N} \left[ 
        \nabla a(x)\cdot \nabla W_{\varepsilon}^{p-1}+a(x) \Delta W_{\varepsilon}^{p-1}
      \right.\\&\qquad\left.  +(\alpha -b)|x|^{\alpha-2}( (\alpha-2+N)W_{\varepsilon}^{p-1} + x\cdot \nabla W_{\varepsilon}^{p-1})
\right] |u|^p dx\\
&\qquad \qquad +\int_{\R^N}c|x|^{\alpha-2}W_{\varepsilon}^{p-1}|u|^{p}+|x|^\beta W_\varepsilon^{p-1}|u|^p dx\\
&\quad = (p-1)\int_{\R^N} (1+|x|^\alpha) W_{\varepsilon}^{p-1} |u|^{p-2}|\nabla |u||^2dx\\
&\qquad  -\frac{1}{p}\int_{\R^N} \left[ 
        (2\alpha-b)|x|^{\alpha-2}x\cdot \nabla W_{\varepsilon}^{p-1}+(1+|x|^\alpha) \Delta W_{\varepsilon}^{p-1}
   \right.\\&\qquad\left. +(\alpha-b)(\alpha-2+N)|x|^{\alpha-2} W_{\varepsilon}^{p-1}
\right] |u|^p dx \\
&\qquad \qquad +\int_{\R^N}c|x|^{\alpha-2}W_{\varepsilon}^{p-1}|u|^{p}+|x|^\beta W_\varepsilon^{p-1}|u|^p dx.
\end{align*}
We first note that

\begin{align*}
 &\nabla W_{\varepsilon}^{p-1}=\gamma W^p_{\varepsilon}|x|^{-\alpha}x\,,\,\,\text{ and }\\
 &\Delta W_{\varepsilon}^{p-1}
  = |x|^{-\alpha}\gamma
        \left( N-\alpha+p(\alpha-2)W_\varepsilon|x|^{2-\alpha} \right)W_\varepsilon^{p},
\end{align*}
with $\gamma = (\alpha-2)(p-1)$, to infer that

\begin{align}\label{estA0}
&\Rp\left(-\int_{\R^N} A_{b,c} u W_{\varepsilon}^{p-1}\overline u|u|^{p-2}dx\right)\\\nonumber
&\quad \geq (p-1)\int_{\R^N} (1+|x|^\alpha) W_{\varepsilon}^{p-1} |u|^{p-2}|\nabla |u||^2dx\\\nonumber
&\qquad  -\frac{1}{p}\int_{\R^N} \left[ 
        (2\alpha-b)\gamma+(1+|x|^\alpha) \gamma |x|^{-\alpha}\left( N-\alpha+p(\alpha-2)|x|^{-\alpha+2}W_\varepsilon \right)\right.\\\nonumber
&\qquad \qquad\left.       +(\alpha-b)(\alpha-2+N)|x|^{\alpha-2} W_{\varepsilon}^{-1}
\right] W_\varepsilon^{p}|u|^p dx\\\nonumber
&\qquad +\int_{\R^N}c|x|^{\alpha-2}W_{\varepsilon}^{p-1}|u|^{p}+|x|^\beta W_\varepsilon^{p-1}|u|^p dx.\\\nonumber
\end{align}

Now, we begin estimating the term $(p-1)\int_{\R^N} (1+|x|^\alpha) W_{\varepsilon}^{p-1} |u|^{p-2}|\nabla |u||^2dx$
by mean of the Hardy inequality \eqref{eq:gen-hardy-weight}. First, we choose $v=W_\varepsilon^{\frac{p-1}{p}}|u|$ and   
$\delta=\alpha$. \\
One has

\begin{align*}
\left( \frac{\alpha-2+N}{p}\right)^2&\int_{\R^N}|x|^{\alpha-2}W_\varepsilon^{p-1}|u|^p dx
    \leq \int_{\R^N}|x|^\alpha W_{\varepsilon}^{p-1}|u|^{p-2}|\nabla |u||^2dx\\   \nonumber 
&\qquad +\frac{\gamma^2}{p^2}\int_{\R^N} |x|^{-\alpha+2}W_\varepsilon^{p+1}|u|^pdx
    \\&\qquad-\frac{2}{p^2}\int_{\R^N} \left(|x|^\alpha \Delta W_\varepsilon^{p-1}
        +\alpha|x|^{\alpha-2}x\cdot \nabla W_\varepsilon^{p-1}\right)|u|^pdx
\end{align*}
that is 

\begin{align*}\label{hardydelta}
&\int_{\R^N}|x|^\alpha W_\varepsilon^{p-1}|u|^{p-2}|\nabla |u||^2dx\geq\\\nonumber
&\quad \int_{\R^N}
\left[\left( \frac{\alpha+N-2}{p} \right)^2|x|^{\alpha-2}W_\varepsilon^{-1}
    -\frac{\gamma^2}{p^2}|x|^{-\alpha+2}W_\varepsilon
\right.\\\nonumber
&\quad \left.
    +\frac{2}{p^2}\gamma(N-\alpha+p(\alpha-2)|x|^{-\alpha+2}W_\varepsilon)
    +\frac{2}{p^2}\gamma \alpha 
\right]W_\varepsilon^p|u|^p dx.
\end{align*}

In the same way we choose $v=W_\varepsilon^{\frac{p-1}{p}}|u|$ and $\delta=0$ to obtain

\begin{align}
&\left( \frac{N-2}{p}\right)^2\int_{\R^N}|x|^{-2}W_{\varepsilon}^{p-1}|u|^p dx 
    \leq \int_{\R^N} W_{\varepsilon}^{p-1}|u|^{p-2}|\nabla |u||^2dx\\\nonumber
&\qquad +\frac{\gamma^2}{p^2}\int_{\R^N} |x|^{-2\alpha+2}W_{\varepsilon}^{p+1}|u|^pdx
        -\frac{2}{p^2}\int_{\R^N} \Delta W_{\varepsilon}^{p-1}|u|^pdx
\end{align}
that is

\begin{align}\label{hardy0}
&\int_{\R^N} W_\varepsilon^{p-1}|u|^{p-2}|\nabla |u||^2dx\geq \int _{\R^N}|x|^{-\alpha}W_\varepsilon^p|u|^p\\
&\left[
\left( \frac{2-N}{p} \right)^2|x|^{\alpha-2}W_\varepsilon^{-1}
   -\frac{\gamma^2}{p^2}|x|^{-\alpha+2}W_\varepsilon
    +\frac{2}{p^2}\gamma\left( N-\alpha+p(\alpha-2)|x|^{-\alpha+2}W_\varepsilon \right)
\right]dx.
\nonumber
\end{align}

Therefore, putting together \eqref{estA0}, \eqref{hardydelta} and \eqref{hardy0}, we get 

\begin{align*}
&\Rp\left(-\int_{\R^N} A_{b,c} u W_{\varepsilon}^{p-1}\overline u|u|^{p-2}dx\right)\geq
        \int_{\R^N} W_\varepsilon^p|u|^p\left[\right.\\
&\qquad -\frac{1}{p}(2\alpha-b)\gamma-\frac{1}{p}\gamma(N-\alpha)
    +(p-1)\left( \frac{2}{p^2}\gamma(N-\alpha)+\frac{2}{p^2}\gamma\alpha \right)\\
&\qquad +\left(-\frac{1}{p}\gamma p(\alpha-2)
    +(p-1)\left(-\frac{\gamma^2}{p^2}+\frac{2}{p^2}\gamma p(\alpha-2)  \right)  \right)|x|^{-\alpha+2}W_\varepsilon\\
&\qquad +\left(
        -\frac{1}{p}(\alpha-b)(\alpha+N-2)+(p-1)\left(\frac{\alpha+N-2}{p}\right)^2
        +c
        \right)|x|^{\alpha-2}W_\varepsilon^{-1}\\
&\qquad +|x|^{-\alpha}\left( -\frac{1}{p}\gamma(N-\alpha)+ (p-1)\frac{2}{p^2}\gamma(N-\alpha) \right)\\
&\qquad +|x|^{-\alpha}\left( -\frac{1}{p}\gamma p(\alpha-2)
        +(p-1)\left(-\frac{\gamma^2}{p^2}+\frac{2}{p^2}\gamma p(\alpha-2)  \right) \right)|x|^{-\alpha+2}W_\varepsilon\\
&\qquad +|x|^{-\alpha} (p-1)\left( \frac{N-2}{p} \right)^2 |x|^{\alpha-2}W_\varepsilon^{-1}\\
&\qquad \left.+|x|^{\beta-\alpha+2}|x|^{\alpha-2}W_\varepsilon^{-1}\right]dx\\
&\quad =
        \int_{\R^N} W_\varepsilon^p|u|^p\left[\right.\\
&\qquad +\frac{\gamma}{p^2}\left( -2\alpha+(N-\alpha)(p-2) \right)+b\frac{\gamma}{p}\\
&\qquad -\frac{\gamma}{p^2}(\alpha-2)  |x|^{-\alpha+2}W_\varepsilon\\
&\qquad +\left(\frac{\alpha-2+N}{p}
        \left( b+\frac{-\alpha+(p-1)(N-2)}{p}\right)+c\right)|x|^{\alpha-2}W_\varepsilon^{-1}\\
&\qquad +|x|^{-\alpha}
        \left(\frac{\gamma}{p^2}(N-\alpha)(p-2) \right)\\
&\qquad -|x|^{-\alpha}\left(\frac{\gamma}{p^2}(\alpha-2) \right)|x|^{-\alpha+2}W_\varepsilon\\
&\qquad +|x|^{-\alpha}{(p-1)\left( \frac{N-2}{p} \right)^2}|x|^{\alpha-2}W_\varepsilon^{-1}\\
&\qquad \left.+|x|^{\beta-\alpha+2}|x|^{\alpha-2}W_\varepsilon^{-1}\right]dx.
\end{align*}

Now, taking into account that $|x|^{-\alpha+2}W_\varepsilon\leq 1$,
$|x|^{\alpha-2}W_\varepsilon^{-1}=1+\varepsilon|x|^{\alpha-2}$
and that $\gamma\leq 0$, we finally have

\begin{align*}
&\Rp\left(-\int_{\R^N} A_{b,c} u W_{\varepsilon}^{p-1}\overline u|u|^{p-2}dx\right)\geq
        \int_{\R^N} W_\varepsilon^p|u|^p\left[ \right.\\
&\qquad +\frac{\gamma}{p^2}\left( -2\alpha+(N-\alpha)(p-2) \right)+b\frac{\gamma}{p}
        -\frac{\gamma}{p^2}(\alpha-2)  
        \\&\qquad+\frac{\alpha-2+N}{p}\left(b+\frac{-\alpha+(p-1)(N-2)}{p}\right)+c\\
&\qquad +|x|^{-\alpha}\left[
        \frac{\gamma}{p^2}(N-\alpha)(p-2) 
        -\frac{\gamma}{p^2}(\alpha-2) 
        +{(p-1)\left( \frac{N-2}{p} \right)^2}\right]\\
&\qquad \left. 
        +|x|^{\alpha-2}\varepsilon\left[
        \frac{\alpha-2+N}{p}
            \left(b+\frac{-\alpha+(p-1)(N-2)}{p}\right)
            +c\right]
        +|x|^{\beta-\alpha+2}\right]dx\\
&\quad =\int_{\R^N} W_\varepsilon^p|u|^p \left[
    -(\alpha-b+\gamma)\frac{\gamma+\alpha+N-2}{p}+(p-1)\left( \frac{\gamma+\alpha+N-2}{p} \right)^2+c\right.\\
&\qquad +|x|^{-\alpha}\left( -\gamma\frac{\gamma+N-2}{p}
        +(p-1)\left( \frac{\gamma+N-2}{p} \right)^2\right)\\
&\qquad 
        +|x|^{\alpha-2}\varepsilon
            \left[
            \frac{\alpha-2+N}{p}\left(b+\frac{-\alpha+(p-1)(N-2)}{p}\right)
            +c\right]
      \\ &\qquad \left. + |x|^{-2}\varepsilon(p-1)\left(\frac{N-2}{p}\right)^2
        +|x|^{\beta-\alpha+2}\right]dx\\
&\quad = \int_{\R^N} W_\varepsilon^p|u|^p \left[
            k_1+|x|^{-\alpha}k_2+|x|^{\alpha-2}\varepsilon k_3+|x|^{-2}\varepsilon k_4+|x|^{\beta-\alpha+2}\right]dx,
\end{align*}
where, $k_1, k_2, k_3, k_4$ are explicit constants which are given by

\begin{align*}
 k_1&=-(\alpha-b+\gamma)\frac{\gamma+\alpha+N-2}{p}
    +(p-1)\left( \frac{\gamma+\alpha+N-2}{p} \right)^2+c
 \\&=\left( \alpha-2+\frac{N}{p} \right)\left( -\alpha+b+\frac{N}{p'} \right)
    +c,
\end{align*}
\[
 k_2= -\gamma\frac{\gamma+N-2}{p} +(p-1)\left( \frac{\gamma+N-2}{p} \right)^2=
    \left( \frac{\alpha}{p'}-2+\frac{N}{p} \right)\left( -\frac{\alpha}{p'}+\frac{N}{p'} \right),
\]
\[
 k_3=\frac{\alpha-2+N}{p}\left(b+\frac{-\alpha+(p-1)(N-2)}{p}\right)+c,
\]
and
\[k_4=(p-1)\left(\frac{N-2}{p}\right)^2.\]

If $\gamma+N-2>0$ which is equivalent to $p<\frac{N-\alpha}{2-\alpha}$ (that is our natural condition) one has $k_2\geq 0$.

Now, it suffices to add and subtract to the right hand side the term  $a_1\int_{\R^N} W_\varepsilon^{p-1}|u|^pdx$, where $a_1\geq 0$ is a suitable constant to be chosen later. First, using H\"older inequality we obtain 
\[
 \int_{\R^N} W_\varepsilon^{p-1}|u|^pdx\leq 
    \left( \int_{\R^N} \left(W_\varepsilon^{p-1}|u|^{p-1}\right)^{p'}dx \right)^{\frac{1}{p'}}
        \left( \int_{\R^N} |u|^p\,dx\right)^{\frac{1}{p}}=\|W_\varepsilon u\|_p^{p-1}\|u\|_p.
\]
Hence, we get

\begin{align*}
&\Rp\left(-\int_{\R^N} A_{b,c} u W_{\varepsilon}^{p-1}\overline u|u|^{p-2}dx\right)\geq\\
&\quad\int_{\R^N}  \left[
            k_1+|x|^{-\alpha}k_2+|x|^{\alpha-2}\varepsilon k_3+|x|^{-2}\varepsilon k_4+|x|^{\beta-\alpha+2}
            +a_1(|x|^{2-\alpha}+\varepsilon)
            \right]W_\varepsilon^p|u|^pdx\\
&\qquad -a_1\|W_\varepsilon u\|_p^{p-1}\|u\|_p. 
\end{align*}

Now, in order to establish (\ref{estim 1}) and (\ref{estim 2})  we analyze the function
\[
 q(x)=k_1+|x|^{-\alpha}k_2+|x|^{\beta-\alpha+2}+|x|^{2-\alpha}a_1+
 \varepsilon\left(|x|^{\alpha-2} k_3+|x|^{-2} k_4+a_1\right).
\]
We first consider the case $k_2>0$ 
 and claim that for every $k>0$ there exists a suitable $a_1\geq 0$ such that $q(x)\geq k$.
Indeed, observe that since $k_4>0$ and $-2<-\alpha$ we have
$\lim_{x\to 0} |x|^{\alpha-2} k_3+|x|^{-2} k_4=+\infty$
and $\lim_{x\to \infty} |x|^{\alpha-2} k_3+|x|^{-2} k_4=0$;
then the function $|x|^{\alpha-2} k_3+|x|^{-2} k_4$ attains a minimum $v_1$.
Taking $a_1>|v_1|$ we have 
\[
 q(x)\geq k_1+|x|^{-\alpha}k_2+|x|^{\beta-\alpha+2}+|x|^{2-\alpha}a_1.
\]
Now set $k\in \R$ and observe that 
$\lim_{x\to 0}k_1+|x|^{-\alpha}k_2+|x|^{\beta-\alpha+2}-k=+\infty$
and $\lim_{x\to \infty}k_1+|x|^{-\alpha}k_2+|x|^{\beta-\alpha+2}-k=+\infty$;
then the function $k_1+|x|^{-\alpha}k_2+|x|^{\beta-\alpha+2}-k$ attains a minimum
$\nu_2$ and, moreover, there exists $r_1$ such that for $|x|\leq r_1$ 
we have $k_1+|x|^{-\alpha}k_2+|x|^{\beta-\alpha+2}-k>0$.
Then we can choose $a_1$ such that $r_1^{2-\alpha}a_1\geq |\nu_2|$ and
 have that $k_1+|x|^{-\alpha}k_2+|x|^{\beta-\alpha+2}+|x|^{2-\alpha}a_1-k>0$
for every $x\in \R^N$.
Therefore taking $a_1\geq \max\{|\nu_1|, \frac{|\nu_2|}{r_1^{2-\alpha}}\}$ we obtain
the desired result.

As regards the case $k_2=0$ we can argue in the same way taking into account that 
$\lim_{x\to 0}k_1+|x|^{-\alpha}k_2+|x|^{\beta-\alpha+2}-k=k_1.$

\qed

As an easy consequence of the above result one has that $W$ is $A_{0,0}-$bounded that is

\[
 \|Wu\|_p\leq k^{-1}\|A_{0,0}\|+C\|u\|,\quad u\in D(A_{0,0})\subset D(W).
\]

Now, we are in position to prove the main result of this section.
\begin{theorem}\label{th:accrev-A_0+W} Assume that $N\geq3$, $0 \le \alpha < 2$, $\beta>0$ and  let $\mathcal{D}=C_c^\infty(\R^N)$.
\begin{itemize}
\item[$\left(i \right)$] If $1<p<\frac{N-\alpha}{2-\alpha}$, then for every $c\in \R$,
$ A_{0,c}$ is quasi $m-$accretive on $ D_p(A_{0,0})$ and  $\mathcal D$ is a core.\\
\item[$\left(ii \right)$] If $p=\frac{N-\alpha}{2-\alpha}$ and if $c>-k_1$,
where 
$k_1=\left( \alpha-2+\frac{N}{p} \right)\left( -\alpha+\frac{N}{p'} \right)$, then 
$ A_{0,c}$ is quasi $m-$accretive on $ D_p(A_{0,0})$ and  $\mathcal D$ is a core.
Moreover $ A_{0,k_1}$ is essentially m-accretive.
\end{itemize}
\end{theorem}
{\sc Proof.}
In order to apply Theorem \ref{Okazawa 1.7}, we observe that both operators $A_{0,0}$ and $W$ are m-accretive in $L^p(\R^N)$. 
Then Lemma \ref{estimate-pot-eps},
with $b=c=0$,
yields $\left(i\right)$ in Theorem \ref{Okazawa 1.7} with $a_2 = 0$,
moreover if $b=c=0$ we have that $k_1>0$.
Observe that  $\left(i\right)$ holds also for $-A_{0,0}+\lambda$ for $\lambda>0$ by changing the value of $a_1$.

The second assumption $\left(ii\right)$ in Theorem \ref{Okazawa 1.7} is obviously satisfied. The last one,
$\left(iii\right)$, follows from  Theorem \ref{thmform}. 

\qed

\subsection{The general case}
In the following subsection we will prove the essential maximality of $(A_{b,c},\mathcal D)$. The first result is the density in $L^{p}(\R^N)$ of the range of $\lambda - A_{b,c}$  for some $\lambda>0.$   
\begin{proposition}\label{pr:maximalityA}
Let $0\le \alpha<2$ and $1<p<\frac{N-\alpha}{2-\alpha}$
or $p=\frac{N-\alpha}{2-\alpha}$ and

$$b\left( \frac{\alpha+N-2}{2} \right)+\left(\frac{N}{p}-2+\alpha  \right)\left( \frac{N}{p'} -\alpha\right)+c>0.$$
Then the range $(\lambda -A_{b,c})\left(\mathcal D\right)$ is dense in $L^{p}(\R^N)$ for some $\lambda>0.$
\end{proposition}
{\sc Proof.}
Suppose that there exists $f\in L^{p\prime}\left( \R^N \right)$ such that
\[
\int_{\R^N} \overline{f}(\lambda-A_{b,c})\varphi \,dx=0\ \text{ for all }\varphi 
	\in \mathcal D.
\]
We have to show that $f = 0$.

Let $\phi=(1+|x|^\alpha)^{b/\alpha}$, where $b\in \R$ is the coefficient of the drift term of $A_{b,c}$
and set  $\varphi=\frac{v}{\sqrt \phi}\in \mathcal D$ with $v\in \mathcal D$.
We note that the function $\phi$ is the function for which the following holds
$$
\frac{1}{\phi}{\rm div}\left(\phi\nabla \varphi\right)=\Delta \varphi +b\frac{|x|^{\alpha -2}}{1+|x|^\alpha}x\cdot \nabla \varphi,\quad \varphi\in \mathcal D.
$$

A simple computation gives

\begin{align}\label{eq:lambdau-Au}
& \left(-A_{b,c}\right)\varphi=\frac{1}{\sqrt \phi}\left[-(1+|x|^\alpha)\Delta v+Uv\right],
\end{align}
where

\begin{align*}
U=(1+|x|^\alpha)\left[-\frac{1}{4}\left|\frac{\nabla \phi}{\phi}\right|^2+\frac{1}{2}\frac{\Delta \phi}{\phi}\right]+{|x|^{\beta}+c|x|^{\alpha-2}},
\end{align*}
that is 

\begin{align*}
&U=
\tilde c|x|^{\alpha-2}-\frac{|x|^{2\alpha-2}}{1+|x|^\alpha}\left( -\frac{b^2}{4}+\frac{b\alpha}{2} \right)+|x|^\beta
\end{align*}
with
$\tilde c=b\left( \frac{\alpha +N-2}{2} \right) +c$.

Now by \eqref{eq:lambdau-Au} we have

\[
\int_{\R^N} \overline{f}(\lambda  -A_{b,c})\varphi\,dx =\int_{\R^N} \lambda\overline{f}\varphi
	+\frac{\overline{f}}{\sqrt \phi}\left[-(1+|x|^\alpha)\Delta v+Uv\right]dx =0.
\]
We define $h=\frac{f}{\sqrt \phi}\in L_{loc}^{p'}(\R^N)\subset L^1_{loc}(\R^N)$, then

\[
\int_{\R^N} 
	\overline{h}\left[-(1+|x|^\alpha)\Delta +U+\lambda \right]v\,dx =0.
\]
Hence

\[
\int_{\R^N} \overline{h}(1+|x|^\alpha)\Delta v\,dx
    =\int_{\R^N}\overline{h}v\left(\lambda+U\right)dx
\]
and since $p(\alpha-2)+N>\alpha>0$ we have that $|x|^{\alpha-2}\in L^p_{loc}(\R^N)$
and then
$\lambda+U\in L^p_{loc}(\R^N)$. Then we have

\[
\Delta \left[h(1+|x|^\alpha)\right]=h\left(\lambda+U\right).
\]
 Applying the Kato's inequality, that is, 

\begin{equation*}\label{K-ineq}
\Delta\vert u\vert \ge \Rp [\sign \overline{u} \Delta u] 
\end{equation*}
in the sense of distributions, for every $u \in L_{loc}^1(\R^N)$ such that $\Delta u  \in L_{loc}^1(\R^N)$,
to $h(1+|x|^\alpha)$ we obtain for every $v\in \mathcal D$

\[
\int_{\R^N} |h(1+|x|^\alpha)|\Delta v\,dx\geq \Rp\int_{\R^N} ({\rm sign}\, \overline{h})
v\Delta \left[h(1+|x|^\alpha)\right]\,dx. \]
Then up to change $\lambda$ for a fixed $\varepsilon>0$ we have

\begin{align*}
&({\rm sign}\, \overline{h})\Delta \left[h(1+|x|^\alpha)\right]
	=|h|\left(\lambda+U\right)\geq |h|\left(\lambda+(\tilde c-\varepsilon)|x|^{\alpha-2}+|x|^\beta\right).
\end{align*}
Then 

\begin{align*}
\int_{\R^N} |h|(1+|x|^\alpha)\Delta v\,dx \geq \int_{\R^N} |h|\left(\lambda+(\tilde c-\varepsilon)|x|^{\alpha-2}+|x|^\beta\right)v\,dx, \end{align*}
and hence 

\begin{align*}
\int_{\R^N} |h|(\lambda-A_{0,\tilde c-\varepsilon})v\,dx  \le 0
\end{align*}
for all $v\in \mathcal D$. 
By Theorem \ref{th:accrev-A_0+W} we have that 
$\overline{\left(\lambda -A_{0,\tilde c-\varepsilon} \right)\mathcal D}=L^p\left( \R^N \right)$ for some $\lambda>0$
if $1<p<\frac{N-\alpha}{2-\alpha}$.
If $p=\frac{N-\alpha}{2-\alpha}$ we require

\[
 b\left( \frac{\alpha+N-2}{2} \right)+\left(\frac{N}{p}-2+\alpha  \right)\left( \frac{N}{p'} -\alpha\right)+c>\varepsilon.
\]

Hence

\[\int_{\R^N} |h|g\,dx  \le 0\]
for all $g\in L^p(\R^N).$ Take $g_n=\sqrt \phi|f|^{p'-1}\vartheta_n\in L^p(\R^N)$ where 
 $\vartheta_n$ is a cutoff function such that $0\leq \vartheta\leq 1$, $\vartheta=1$ in $B(n)$ and $\vartheta=0$ in $B(2n)^c.$ It holds
 
 \[
  \int_{\R^N} |f|^{p'}\vartheta_n\,dx\leq 0
 \]
for all $n\in \N$. This implies that $f=0$ and then the maximality is proved.

\qed                   

Now, we are ready to prove the main result of this section. 
\begin{theorem}\label{generation}
Let $0\leq \alpha<2$ and either $1<p<\frac{N-\alpha}{2-\alpha}$
or else $p=\frac{N-\alpha}{2-\alpha}$ and
\begin{equation}\label{eq:condition-bc}
 \left(\frac{N}{p}-2+\alpha  \right)\left( \frac{N}{p'} -\alpha+b\right)+c>0.
\end{equation}
Then the closure of $(A_{b,c},\mathcal D)$
generates an analytic $C_0-$semigroup.
\end{theorem}
{\sc Proof.}
We observe that since  $1<p\leq \frac{N-\alpha}{2-\alpha}$ we have that $p(\alpha-1)>-N$ and $p(\alpha-2)>-N$,
then the drift and the potential terms of $A_{b,c}$ belong to $L^p(\R^N)$. 
As a consequence  $A_{b,c}$ maps $\mathcal D$ into $L^p(\R^N)$ and $A_{b,c}$ is well defined in $\mathcal D$ as an operator acting on $L^p(\R^N).$

First we consider the case $1<p<\frac{N-\alpha}{2-\alpha}$.
 By Theorem \ref{thmform}, for $\lambda$ greater than a suitable $\lambda_0>0$, the operator $-A_{b,c}+\lambda$ is accretive, then it is closable,
by Proposition \ref{pr:maximalityA}, the range $(-A_{b,c}+\lambda)({\mathcal D})$ (we can take $\lambda>\lambda_0$)
is dense in $L^p(\mathbb{R}^N)$. 
Hence the closure of $(A_{b,c},{\mathcal D})$ is quasi $m-$accretive. By Lumer and Phillips Theorem, cf. e.g. \cite[Ch. I, Theorem. 4.3]{pazy2} or \cite[Theorem 1.1]{Eberle}, the closure of $A_{b,c}$ is the generator of a $C_0-$semigroup on $L^p(\R^N)$.    
Since  by   Theorem \ref{formgeneration}  $A_{b,c}$ is sectorial, then the semigroup is analytic.

Now consider the case $p=\frac{N-\alpha}{2-\alpha}$. 
Fix $b\in \R$ and choose $c_1$ such that the hypotheses of Proposition \ref{pr:maximalityA} are satisfied and such that the constant $k_1(b,c_1)$ of 
Lemma  \ref{estimate-pot-eps}
is positive. Arguing as before we have that
the closure of  $(-A_{b,c_1},{\mathcal D})$ 
is quasi $m-$accretive.
We apply once again Okazawa's perturbation Theorem 
\ref{Okazawa 1.7}
to the operator obtained by the closure of  $(-A_{b,c_1},{\mathcal D})$
perturbed by $W$.
Since $k_1(b,c_1)>0$, Lemma  \ref{estimate-pot-eps} gives condition $(i)$ of Theorem \ref{Okazawa 1.7}.
The second assumption $\left(ii\right)$ is obviously satisfied. The last one,
$\left(iii\right)$, is obtained by  Theorem \ref{thmform}.
Finally we obtain that $-{A_{b,c_1}}+c_2W$ defined on $\mathcal D$ is quasi essentially m-accretive if $c_2\geq -k_1(b,c_1)$.

Setting $c=c_1+c_2$ we have that $c_2\geq -k_1(b,c_1)$ is equivalent
to $k_1(b,c)>0$ which is the condition \eqref{eq:condition-bc}.

\qed

\begin{remark}
In the case $p=\frac{N-\alpha}{2-\alpha}$ we observe that the condition 

$$\left(\frac{N}{p}-2+\alpha\right)\left( \frac{N}{p'} -\alpha+b\right)+c>0,$$
is the same condition as  G. Metafune et al. in \cite[Theorem 4.5]{met-oka-sob-spi}.
\end{remark}

We observe that since $\left(A_{b,c},\mathcal D_0\right)$ is contained in $(A_{b,c},{\mathcal D})$ we have that the closure of $(A_{b,c},{\mathcal D})$ coincides with $ (A_{b,c},D_p(A_{b,c})).$

Finally, in the end of this section, we propose  to prove that the closure of $(A_{b,c},{\mathcal D}_0)$ coincides with $(A_{b,c}, D_p(A_{b,c}))$ and to give more information on the domain $D_p(A_{b,c})$. 
Observe that, by Lemma  \ref{estimate-pot-eps}  and Theorem \ref{Okazawa 1.7}, for every  $u\in  D_p(A_{b,c})$ we have that    
 
\begin{equation}\label{eq:stima-pot-W}
 \|Wu\|_p\leq C\left(\|A_{b,c}u\|_p+\|u\|_p\right)
\end{equation}
for some positive constants $C$.\\

 The following proposition proves that ${\mathcal D}_0$ is a core for $\overline{(A_{b,c},{\mathcal D})}$ and $ D_p(A_{b,c})$ is included in some weighted Sobolev space.       

\begin{proposition}
If the hypotheses of Theorem \ref{generation} hold, then 
${\mathcal D}_0$ is a core for $(A_{b,c}, D_p(A_{b,c}))$ which coincides with the closure of $(A_{b,c},{\mathcal D})$.
Moreover

\begin{align*}
 D_p(A_{b,c})  \subset  \{u\in L^p(\R^N)\cap W^{2,p}_{loc}(\R^N\setminus\{0\}):\,\,
 (1+|x|^\alpha)D^2u, | |x|^{\alpha-1}\nabla u|, |x|^\beta u, Wu  \in L^p(\R^N)\}.
\end{align*}
\end{proposition}
{\sc Proof.}
First we prove that there exists a constant $C>0$ such that for every $u\in {\mathcal D_0}$ the following estimates hold

\begin{align}
 \|(1+|x|^\alpha)\Delta u\|_p\leq C\left( \|A_{b,c}u\|_p+\|u\|_p \right);\label{eq:stima-diff} \\
 \||x|^{\alpha-1}\nabla u\|_p\leq C\left( \|A_{b,c}u\|_p+\|u\|_p \right); \label{eq:stima-drift}\\
 \||x|^{\alpha-2}u\|_p\leq C\left( \|A_{b,c}u\|_p+\|u\|_p \right); \nonumber\\
 \||x|^\beta u\|_p\leq C\left( \|A_{b,c}u\|_p+\|u\|_p \right).\nonumber
\end{align}
Indeed ${\mathcal D}_0\subset  D_p(A_{b,c})$, then if $u\in {\mathcal D}_0$ we have by \eqref{eq:stima-pot-W} and by \cite[Corollary 2.4]{lor-rhan}   
 \begin{equation}\label{eq:stima-pot-VW}
\left\{
\begin{array}{ll}
 \|Wu\|_p\leq C\left( \|A_{b,c}u\|_p+\|u\|_p \right), \\ \ \\
  \||x|^\beta u\|_p\leq C \left( \|A_{0,0}u\|_p+\|u\|_p \right),
\end{array}\right.
\end{equation}
for some positive constant $C$. 
Moreover by \cite[Lemma 4.4]{Met-Sob-Spi} for every $\varepsilon>0$ there exists $C_\varepsilon>0$ such that

\begin{equation}\label{eq:interp-weight}
 \|| x|^{\alpha-1} \nabla u\|_p\leq \varepsilon \| |x|^\alpha D^2u \|_p+C_\varepsilon \|x^{\alpha-2}u\|_p.
\end{equation}
Then applying the Calder\'on-Zygmund inequality to the function $|x|^\alpha u$ and taking into account \eqref{eq:stima-pot-VW} 
we have

\begin{align*}
 &\||x|^\alpha D^2u \|_p\leq C\left( \|D^2(|x|^\alpha u)\|_p+\||x|^{\alpha-1}\nabla u\|_p+\||x|^{\alpha-2}u\|_p \right)\\
 &\quad \leq C\left( \||x|^{\alpha}\Delta u\|_p+\||x|^{\alpha-1}\nabla u\|_p+\||x|^{\alpha-2}u\|_p \right)\\
 &\quad \leq C\left( \|A_{b,c}u\|_p+\||x|^{\alpha-1}\nabla u\|_p+\||x|^\beta u\|_p+\|Wu\|_p \right)\\
 &\quad \leq C\left( \|A_{b,c}u\|_p+\||x|^{\alpha-1}\nabla u\|_p+\|A_{0,0}u\|_p+\|A_{b,c}u\|_p \right)\\
 &\quad \leq C\left( \|A_{b,c}u\|_p+\||x|^{\alpha-1}\nabla u\|_p\right).
\end{align*}
Then by \eqref{eq:interp-weight} 

\begin{align*}
 \| |x|^{\alpha-1} \nabla u\|_p\leq \varepsilon \left( \|A_{b,c}u\|_p+\||x|^{\alpha-1}\nabla u\|_p\right)+C_\varepsilon \|Wu\|_p
\end{align*}
for which the estimate \eqref{eq:stima-drift} and then \eqref{eq:stima-diff} follow.

Let $u\in {\mathcal D}$ and set $u_n=\varphi_n u\in {\mathcal D}_0$,
where $\varphi_n\in \mathcal D_0$ such that $0\leq \varphi_n\leq 1$, $\varphi_n=0$ in $B\left( \frac{1}{n} \right)$ and $\varphi_n=1$ in $B^c\left( \frac{2}{n} \right)$. Moreover $|\nabla \varphi_n(x)|\leq C \frac{1}{|x|}$ and
$|D^2\varphi_n (x)|\leq C \frac{1}{|x|^2}$.
Arguing as in Proposition \ref{pr:core-for-a} we have that $u_n$ converges to $u$ with respect to the operator norm, and then
${\mathcal D}_0$ is a core for the closure of $(A_{b,c},{\mathcal D})$.

Finally we observe that
by local elliptic regularity  it follows that $ D_p(A_{b,c})\subset W_{loc}^{2,p}(\R^N\setminus\{0\})$.
And estimates \eqref{eq:stima-diff}, \eqref{eq:stima-drift} and \eqref{eq:stima-pot-VW} hold in $ D_p(A_{b,c})$ by density.
\qed

\bibliographystyle{amsplain}
\bibliography{bibfile}

\end{document}